\documentclass[draft,12pt]{amsproc}
\usepackage{amscd,amsthm,amsfonts,amssymb,latexsym}
\usepackage{mathrsfs}
\theoremstyle{plain}
\newtheorem*{theorem}{Theorem}
\newtheorem*{cor}{Corollary}

\newtheorem{prop}{Proposition}

\theoremstyle{definition}
\newtheorem{defn}{Definition}

\theoremstyle{remark}

\def\balpha{\mbox{\boldmath$\alpha$}}
\def\bbeta{\mbox{\boldmath$\beta$}}
\def\bgamma{\mbox{\boldmath$\gamma$}}
\def\bdelta{\mbox{\boldmath$\delta$}}
\def\bzeta{\mbox{\boldmath$\zeta$}}
\def\bita{\mbox{\boldmath$\eta$}}
\def\btheta{\mbox{\boldmath$\vartheta$}}
\def\bkappa{\mbox{\boldmath$\kappa$}}
\def\bnu{\mbox{\boldmath$\nu$}}
\def\bxi{\mbox{\boldmath$\xi$}}
\def\bpi{\mbox{\boldmath$\pi$}}
\def\brho{\mbox{\boldmath$\varrho$}}
\def\bsigma{\mbox{\boldmath$\sigma$}}
\def\btau{\mbox{\boldmath$\tau$}}

\def\ma{\mathbb}
\def\mc{\mathcal}
\def\ms{\mathscr}

\def\C{{\ma C}}
\def\N{{\ma N}}
\def\Z{{\ma Z}}
\def\res{\mathop{\rm Res\,}\limits}
\def\sym{\mathop{\rm Sym}\nolimits}

\begin{document}
\title[Ideal Numbers as maximal Ideals]
{Pr\"ufer's Ideal Numbers as Gelfand's maximal Ideals}

\author{S.~Albeverio}

\address{\scriptsize Institut f\"{u}r Angewandte Mathematik,
Universit\"{a}t Bonn, Wegelerstr. 6, D-53115 Bonn, Germany; SFB
611, Bonn, BiBoS, Bielefeld -- Bonn}

\email{albeverio@uni-bonn.de}
\author{V.~Polischook}

\address{\scriptsize\noindent Dept. of Mathematics, St.~Petersburg
State Polytechnical University, Polytechnicheskaya 29, 195251,
St.~Petersburg, Russia.}

\email{polischook@list.ru}

\thanks{This paper was supported in part by DFG Project 436 RUS
113/809/0-1. The second author (V.~P.) was also supported in part
by Grant 05-01-04002-NNIOa of Russian Foundation for Basic
Research.}

\subjclass[2000]{Primary: 43A60 Almost periodic
functions on groups, Secondary: 43A75 Analysis on specific compact
groups}

\date{01.03.2007}

\keywords{Almost periodic number-theoretical functions,
compactification of, Banach algebra, Gelfand's theory,
maximal ideals, Novoselov's theory, polyadic numbers,
polyadic analysis}

\begin{abstract}
Polyadic arithmetics is a branch of mathematics related to
$p$--adic theory. The aim of the present paper is to show that
there are very close relations between polyadic arithmetics and
the classic theory of commutative Banach algebras. Namely, let
$\ms A$ be the algebra consisting of all complex periodic
functions on $\Z$ with the uniform norm. Then the polyadic
topological ring can be defined as the ring of all characters
$\ms A\to\C$ with convolution operations and the Gelfand topology.
\end{abstract}
\maketitle
\section*{Prolegomena}
Let us start from an elementary polynomial identity. Let $\{p_j\}$
be an arbitrary sequence of natural numbers and $z$ be a complex
variable. It is easily seen that for each $k\in\N$
$$
(1+z+\cdots+z^{p_1-1})(1+z^{p_1}+\cdots+z^{p_1(p_2-1)})
(1+z^{p_1p_2}+\cdots+z^{p_1p_2(p_3-1)})\cdots=
$$
$$
=1+z+z^2+z^3+\cdots+z^{p_1p_2\cdots p_k-1}\,.
$$
In particular, if the sequence $p_j$ is constant, i.e. $p_j=p$
for all $j$, the identity will be transformed to the form
\begin{equation*}
\prod_{r=0}^{k-1}\Bigl(\sum_{\nu=0}^{p-1}z^{\nu p^r}\Bigr)=
\sum_{n=0}^{p^k-1}z^n
\end{equation*}
showing that each non--negative integer $n<p^k$ can be uniquely
represented by the sum
\begin{equation}\label{p-adic}
n=\nu_0(n)\cdot p^0 +\nu_1(n)\cdot p^1+\cdots+\nu_{k-1}(n)\cdot
p^{k-1}\,,
\end{equation}
where \emph{the digits} $\nu_r(n)$ satisfy the following
condition: $0\le\nu_r(n)<p$.
\par\noindent
If we set $p_j=j$ then $p_1\cdots p_r=r!$ and the identity has
the form
$$
\prod_{r=0}^{k-1}\frac{1-z^{(r+1)!}}{1-z^{r!}}=
\prod_{r=0}^{k-1}\Bigl(\sum_{\nu=0}^rz^{\nu r!}\Bigr)=
\sum_{n=0}^{k!-1}z^n
$$
which shows that each non--negative integer $n<k!$ admits the
unique \emph{factorial representation}
\begin{equation}\label{factorial}
n=\nu_1(n)\cdot 1!+\nu_2(n)\cdot 2!+\cdots+\nu_{k-1}(n)\cdot(k-1)!
\end{equation}
where \emph{the digits} $\nu_r(n)$ satisfy the condition
$0\le\nu_r(n)\le r$.
\par
Trying to give a general meaning to formulae of the type
(\ref{p-adic}) and (\ref{factorial}) with an \emph{infinite}
number of digits caused H.\,Pr\"ufer \cite{Pruefer} to create a
new theory of \emph{ideal numbers}. J.~von Neumann \cite{Neumann}
further developed this theory which was later completed in the
theory of adeles (\cite{Gelf-and}).
\par
The representation (\ref{p-adic}) of ideal numbers with the prime
$p$ as the radix and an infinite number of digits (invented by
K.~Hensel~\cite{Hensel}) gave rise to the $p$--\emph{adic
analysis} (for references see e.g.~\cite{Robert}).
\par
The factorial representation (\ref{factorial}) with an infinite
number of digits lays at the foundation of the \emph{polyadic
analysis}. M.\,D.~Van Dantzig \cite{VanDantzig} and later
E.~Novoselov \cite{Nov-66} gave an essential contribution to
this area. The most significant results were summarized by
E.~Novoselov in his book \cite{Nov-82}.
\par
By a polyadic number E.~Novoselov means a \emph{finite integral
adele}. Such a number can be viewed as a limit of an integer
sequence which converges $p$--adically for every prime $p$. The
Novoselov's theory includes the polyadic topology, arithmetics, the
basic analytic functions, measure theory and integration in the
polyadic domain.
\par
J.--L.~Mauclaire well adapted and completed this theory in a
number of papers
\cite{Mauclaire-81,Mauclaire-85,Mauclaire-86,Mauclaire-88}.
\par
The most important property of the polyadic topological ring $\ms
P$ with the divisibility topology $\bsigma$ is to be a
\emph{compact Hausdorff space}. As is well known, there is no
translation--invariant $\sigma$--additive measure on $\Z$. But
such a measure (Haar measure) exists on the compact additive group
of the polyadic ring $(\ms P,\bsigma)$ including $\Z$ as dense
subset. So, the \emph{compactification} of the ring $\Z$ is the
key idea of the Novoselov theory. His method for the
investigation of the distribution problems of arithmetic functions
has been studied by several papers and books, such as\,
\cite{Knopfmacher,Krizhyus, Schwarz,Kubota,Sugita}.
\par
Of course, various ways of compactification are possible.
They lead us to different (a priori) completions of $\Z$.
\par
Z.~Krizhyus \cite{Krizhyus} deals with the compactification based
on \emph{periodicity}, which seems very natural. Indeed, on the
one hand, the sum (or the product) of $m$-- and $n$--periodic
functions on $\Z$ has each common multiple of $m$ and $n$ as
period. On the other hand, if a function $f$ is $m$-- and
$n$--periodic simultaneously then any common divisor of $m$ and
$n$ is also a period of $f$. As divisibility is a main notion of
arithmetics, the algebra $\ms A$ of all periodic function on
$\Z$ must play an important role in number theory.
\par
In fact, there are many known works (cf. \cite{DeBruijn,Delsarte,
Knopfmacher,Daboussi,Schwarz,Kubota,Sugita}) concerned with the
\emph{almost periodic} functions which are elements of closures
$\ms A$ with respect to various norms. In \cite{Krizhyus} the
\emph{uniform} completion $\mc A$ of the algebra $\ms A$ is
considered. By the theory of commutative Banach algebras
presented, e.\,g., in \cite{Gelfand}, $\mc A$ is isometrically
isomorphic to the algebra of all continuous functions on the
\emph{Gelfand $\mc A$--compactification $(\ms G,\bgamma)$} of
$\Z$.
\par

The plan of this article is as follows. In the sections 1--3 we
introduce the basic concepts and notations used in the paper.
Then, in the section~4 we define the ring operations in $\ms G$.
Further, in the section~5 we define a \emph{cluster topology} in
the ring $\ms G$ and we prove that it is identical with the
Gelfand topology~$\bgamma$. Next, in the sections 6--7 the
polyadic topological ring $(\ms P,\bsigma)$ is described. This
exposition differs from Novoselov's one merely in non--essential
technical details. The contents of the section 8 can be regarded
as a presentation of some ``bridge'' between the topological rings
$(\ms P,\bsigma)$ and $(\ms G,\bgamma)$. At last, in the section~9
we state the identity of these topological rings.

\section {Preliminaries}

\subsection {Weak topologies in the dual space}

Let $B_0$ be a dense subspace of a Banach space $B$. For any
$\psi\in B^*$, $u\in B_0$ and
$R>0$ we define a \emph{neighborhood $U_R(u;\psi)$} of the element
$\psi$:
$$
U_R(u;\psi)=\{\phi\in B^*:\,|\phi(u)-\psi(u)|<R\}\,.
$$
These neighborhoods form a subbase of some topology in $B^*$
called \emph{$B_0$--topology.} If $B_0=B$ then it
is called the \emph{$*$--weak topology in $B^*$}.

\begin{prop}\label{subbasis}
All topologies induced by different $B_0$--topologies on
a bounded subset $\ms F\subseteq B^*$ coincide.
\end {prop}

\begin{proof}
Let us check that a topology $\btau$ induced on $\ms F$ by an
arbitrary $B_0$--topology coincides with the topology $\bsigma$
induced on this set by the $*$--weak topology of the space $B^*$.
Since it is obvious (from the definition of the $B_0$--topology)
that $\btau\subseteq\bsigma$, it is necessary to show that
$\bsigma\subseteq\btau$. It will suffice to prove that for each
functional $\psi\in\ms F$ an arbitrary element $U_R(u;\psi)$ of
the $\bsigma$--subbase is wider than some element $U_r(v;\psi)$
from the $\btau$--subbase. Let the norms of functionals from $\ms
F$ be bounded by $N$. Let us choose an element $v$ of the set
$B_0$ (dense in $B$) and a positive number $r$ such that the
following inequalities
$$
3\|u-v\|<R/N, \qquad 3r<R
$$
hold. Under this condition
$U_r(v;\psi)\subseteq U_R(u;\psi)$. Indeed, for any functional
$\phi$ from $U_r(v;\psi)$ we have
\begin{multline*}
|\phi(u)-\psi(u)|\le|\phi(u-v)|+|\phi(v)-\psi(v)|+|\psi(v-u)|\le\\
\le\|\phi\|\cdot\|u-v\|+|\phi(v)-\psi(v)|+\|\psi\|\cdot\|u-v\|<\\
<N\cdot\|u-v\|+r+N\cdot\|u-v\|<R/3+R/3+R/3\,,
\end{multline*}
i.\,e. every functional $\phi$ from $U_r(v;\psi)$ belongs to
$U_R(u;\psi)$.
\end{proof}

\subsection {Topological ring} \label{base}

A subset $U$ of a commutative ring $R$ is called \emph{symmetric}
if with each element $a\in U$ it contains its opposite $(-a)$. In
order to define in $R$ a Hausdorff topology compatible with the
algebraic structure of $R$, it is convenient to set a \emph{base
of neighborhoods of zero} that is a class $\ms B$ of symmetric
subsets complying with the following conditions:

\begin{enumerate}
\item The zero of the ring $R$ is a unique common element of all
      neighborhoods.
\item An intersection of two arbitrary neighborhoods of zero
      includes some neighborhood of zero.
\item For each neighborhood of zero $V$ it is possible to
      find a neighborhood of zero $U$ such that
      \quad $\{u_1+u_2:u_1,u_2\in U\}\subseteq V$.\label{U+U}
\item For each neighborhood of zero $V$ and arbitrary element
      $a\in R$ it is possible to show a neighborhood $U$ such that
      $\{au:u\in U\}\subseteq V$.
\item For an arbitrary neighborhood of zero $V$ it is possible
      to specify a neighborhood of zero $U$ such that
      \quad $\{u_1u_2:u_1,u_2\in U\}\subseteq V$.
\end{enumerate}

The \emph{translations} $U_a=\{a+u:u\in U\}$ of neighborhoods of
zero form a base $\ms B_a$ of neighborhoods of an arbitrary point
$a$:
$$
\ms B_a=\{U_a:U\in\ms B\}\,.
$$
An empty set and any subset $S\subseteq R$ including some
neighborhood $U_a$ of each its point $a$ are announced as
\emph{open} in the topology of the ring $R$. Obviously, by this
definition, a union of open sets and an intersection of two open
sets are open.

\begin{defn}
\emph{The interior} $W_a^o$ of a neighborhood $W_a$ consists of
points $U$ having a neighborhood of zero $V$ such that
$V_u\subseteq W_a$.
\end{defn}
\noindent
The interior $W_a^o$ is not empty, because it contains the
point $a$.

\begin {prop}
The set $W_a^o$ is the largest open part of $W_a$.
\end {prop}

\begin {proof}
Let $u\in W_a^o$ and let $V$ be a neighborhood of zero for which
$V_u\subseteq W_a$. Let $U$ be a neighborhood of zero found
according to (\ref{U+U}). Let us check that $U_u\subseteq W_a^o$
and thus establish that $W_a^o$ is open.
\par
The neighborhood $U_v$ of an arbitrary point $v\in U_u$
consists of points $u+u_1+u_2$ where $u_1$, $u_2$
belong to the neighborhood $U$. By virtue of (\ref{U+U}),
$u_1+u_2\in V$, so that $U_v\subseteq V_u\subseteq W_a$.
Hence, the arbitrary point $v$ of the neighborhood $U_u$
belongs to $W_a^o$, as was to be checked.
\par
Let now $S$ be a nonempty open subset of a neighborhood $W_a$.
Each point $u$ belongs to $S$ together with some
neighborhood $V_u$. All the more, $u$ belongs to $W_a$
together with $V_u$, i.\,e. $u\in W_a^o$.
\end {proof}

\begin {prop}
Let topologies $\bsigma$ and $\btau$ in $R$ be defined by two bases
of neighborhoods $\ms A$ and $\ms B$, respectively.
The topology $\bsigma$ is weaker than the topology $\btau$
if and only if each neighborhood $U\in\ms A$ includes some
neighborhood $V\in\ms B$.
\end {prop}

\begin{proof}
Let $\bsigma$ be weaker than $\btau$, i.\,e. every $\bsigma$--open
set is $\btau$--open. The interior $U^o$ of the neighborhood $U$
contains the zero of the ring $R$ and it is $\bsigma$--open.
A fortiori it is $\btau$--open, so that zero is included into it
with some $\btau$--neighborhood $V$. Thus,
$V\subseteq U^o\subseteq U$.
\par
Conversely, each point $s$ of the $\bsigma$--open set $S$ belongs
to it with a neighborhood $U_s$ which is a translation of a
neighborhood of zero $U$. The neighborhood $U$ contains some
neighborhood $V\in\ms B$. Since the above translation keeps the
inclusion of sets $V_s\subseteq U_s\subseteq S$, the point $s$ of
$S$ belongs to $S$ with a $\btau$--neighborhood $V_s$. Hence this
set is $\btau$--open.

\end{proof}
\section {Almost periodic functions on $\Z$}

\subsection {Algebra $\mc A$ of almost periodic functions}
\label {almost}
Denote by $\ms A$ the class of all complex--valued periodic
functions $u$ on $\Z$. This is a commutative algebra with the unit
$e$, termwise operations and a complex conjugation $u\mapsto{u^*}$
as involution. The class $\ms A_p$ consisting of all $p$--periodic
complex--valued functions on $\Z$ is a finite--dimensional
subalgebra in $\ms A$. Let $p$ and $q$ be arbitrary positive
integers. Any $q$--periodic function is also a $pq$--periodic one.
Hence $\ms A$ can be considered as an \emph{inductive limit} of
the increasing sequence of algebras
$$\ms A_{1!}\subseteq\ms A_{2!}\subseteq\cdots
\subseteq\ms A_{n!}\subseteq\cdots\,.
$$
A periodic function $u$ on $\Z$ takes only a finite number of
values, hence one can define on $\ms A$ a uniform norm
$\|u\|=\sup|u(n)|$ having the following properties:
$\|e\|=1;\quad\|uv\|\le\|u\|\cdot\|v\|;\quad\|u^*u\|=\|u\|^2$. By
completion of $\ms A$ we get the Banach algebra $\mc A$ of
\emph{almost periodic} functions on $\Z$.

\subsection{Subalgebra of $p$--periodic functions}

\begin{defn}
Denote by $\res_{\!p}(n)$ a \emph{residuum} under division
$n\in\Z$ by $p\in\N$ defined by the equality $n=mp+\res_p(n)$
\, $(0\le n-mp<p)$. We define a \emph{raster} $R_p(m)$ as
the set of level $m$ of the function $\res_p$:

\begin{equation*}
R_p(m)=\{k\in\Z:\res_p(k)=m\}\qquad(0\le m< p)
\end{equation*}
\end{defn}

\begin{prop}\label{indicators}
The indicators $e_m^p$ $(0\le m<p)$ of the rasters $R_p(m)$
have the following properties:
\begin {enumerate}
\item $e_m^pe_n^p=\delta_{mn}e_m^p$;
\item $\sum e_m^p=e$;
\item $u=\sum u(m)e_m^p$ for any $u$ from $\ms A_p$;
\item the functions $e_m^p$ form a base in $\ms A_p$.
\end {enumerate}
\end {prop}

\begin{proof}
\par\noindent
(1) The product of indicators is the indicator of the intersection,
but the sets of different levels do not intersect, therefore
if $m\ne n$ then $e_m^pe_n^p=0$. Moreover by definition any
indicator coincides with its square.
\par\noindent
(2) The rasters $R_p(m)$ (as sets of all different levels) split
$\Z$ into components. But in case of partition, the indicator $e$
of a union is equal to the sum of the component indicators $e_m^p$.
\par\noindent
(3) On the raster $R_p(m)$ the product $(u-u(m)e)\,e_m^p$ is equal
to zero because of the first factor, on its completion it is equal
to zero because of the second one. Therefore $ue_m^p=u(m)e_m^p$
and hence, in view of equality (2),
$$
u=ue=u\sum e_m^p=\sum ue_m^p=\sum u(m)e_m^p\,.
$$
(4) The condition $e_m^p(k)=1$ is equivalent to
$$
\res_p(k)=\res_p(k+p)=m\,.
$$
The latter equality means that $(k+p)\in R_p(m)$.
Thus, $e_m^p(k+p)=1$, and the $p$--periodicity of the indicator
$e_m^p$ is established.
Let us check the independence of the elements $e_m^p$ whose
completeness was proved in the previous item. Multiplying both
parts of the equality $\sum\alpha_me_m^p=0$ by $e_n^p$ we get, in
view of the item $(1)$, the equality $\alpha_ne_n^p=0$ which is
possible only if $\alpha_n=0$.
\end{proof}

\begin{prop}\label{separ}
The subalgebra $\ms A_p$ separates the points of the set
$$
K=\{0,1,\ldots,p-1\}\,.
$$
\end{prop}

\begin{proof}
Let $0\le m<p$, then $e^p_m(k)=1$ for $k=m$, therefore $e^p_m(k)=0$
for all other values $k\in K$.
\end{proof}
\section{Characters of the Banach algebra $\mc A$}

\subsection{General information}

We call by \emph{character} any non--zero multiplicative linear
functional $\psi:\mc A\to\C$ . Such a functional is always bounded
and its norm is equal to $1$. Denote the set of all characters by
$\ms G$. A non--trivial ideal of the algebra $\mc A$ is called
\emph{maximal} if it is not a subset of any other non--trivial
ideal. The kernel of a character is a maximal ideal and
conversely: any maximal ideal is the kernel of a unique character.
\par
For each element $u\in\mc A$ we define a function $\hat u:\ms
G\to\C$ by the formula $\hat u(\psi)=\psi(u)$. The topology
$\bgamma$ induced on $\ms G$ by the $*$--weak topology of the dual
$\mc A^*$ is called \emph{Gelfand topology}. It is the weakest
topology for which all functions $\hat u$ are continuous. Since the
subalgebra $\ms A$ is dense in $\mc A$ and $\ms G$ is a subset of
the unit sphere by virtue of Proposition\,\ref{subbasis}, it is
possible to define the topology $\bgamma$ by such a subbase of
neighborhoods of points $\psi_0\in\ms G$:
\begin{equation}\label{neighborhood}
    U_R(u;\psi_0)=\{\psi\in\ms G:\,|\psi(u)-\psi_0(u)|<R\}
    \qquad (u\in\ms A,\, R>0)\,.
\end{equation}
Owing to classical inferences of the general Gelfand theory,
the topological space $(\ms G,\bgamma)$ is a compact Hausdorff
one and the transform $G:u\mapsto\hat u$ is an isometry of the
algebra $\mc A$ onto the algebra of all continuous complex--valued
functions $C(\ms G,\bgamma)$ with a uniform norm.
\par
In view of the density of $\ms A$ in the Banach algebra $\mc A$
both algebras have the same dual and the same set of characters.
Thus it is sufficient to describe only the characters of the algebra
$\ms A$.

\subsection{Restriction of character on the subalgebra
$\ms A_p$\,}

For all $k\in\N$ we define a functional $\bdelta_k$ by the equality
$\bdelta_k(u)=u(k)$. Its linearity and multiplicativity are obvious,
while its non--triviality follows from the equality
$\omega_k(e)=e(k)=1$.

Denote by $\psi_p$ the restriction $\psi|\ms A_p$ of a character
$\psi$ on the subalgebra $\ms A_p$. As the following proposition
shows, the set of non--zero restrictions $\psi_p$ consists only of
the functionals $\bdelta_k$.

\begin{prop}\label{delta}
For each $p\in\N$ and $\psi\in\ms G$ either the restriction
$\psi_p$ is trivial or there is a unique non--negative number
$k=\bkappa(p;\psi)<p$ such that $\psi_p=\bdelta_k$.
\end{prop}

\begin{proof} For a character $\psi$ with a non--zero
restriction $\psi_p$ there is an element $u$ in the subalgebra
$\ms A_p$ such that $\psi(u)\ne0$. Dividing by $\psi(u)$
both members of the equality $\psi(u)=\psi(eu)=\psi(e)\psi(u)$
we conclude that $\psi(e)=1$.
Further, it follows from the item (1) of
Proposition\,\ref{indicators} that $\psi^2(e_k^p)=\psi(e_k^p)$,
i\,.e. for each value of the index $k$ either $\psi(e_k^p)=0$, or
$\psi(e_k^p)=1$. We find by calculating the value of $\psi$ on each
side of the equality (2) from Proposition\,\ref{indicators} that
$\sum\psi(e_m^p)=1$. This is only possible if exactly one of the
addends differs from zero. So, for the considered character $\psi$
the set $\{k:0\le k<p\}$ contains a unique $k$ such that the
equality $\psi(e_m^p)=\delta_{km}$ holds for any $m$. Let us
compute the value of the functional $\psi$ on an arbitrary
function $u\in\ms A_p$ using the equality (3) from
Proposition\,\ref{indicators}:
$$
\psi(u)=\sum u(m)\psi(e_m^p)=\sum u(m)\delta_{km}=u(k)
=\bdelta_k(u)\,.
$$
Assume now that $\psi=\bdelta_r$ where $0\le r<p$. Then
$u(k)=u(r)$ for every $p$--periodic function $u$. In particular,
for $u=\res_p$ we have\\
\phantom{MMMMMMMMMMMMM} $r=u(r)=u(k)=k$.\hfill
\end {proof}
\noindent
Let us underline an important property of the function $\bkappa$.

\begin{prop}\label{WOW!}
If the restrictions of a character $\psi$ to the subalgebras
$\ms A_m$ and $\ms A_n$ are non--trivial, then
$$
\bkappa(m;\psi)\equiv\bkappa(n;\psi)\mod\ (m\land n)
$$
where $m\land n$ is the greatest common divisor of $m$ and $n$.
\end{prop}

\begin{proof}
Set $m\land n=d$. As the function $\res_d$ belongs to both
$\ms A_m$ and $\ms A_n$, we have
$$
\res_d(\bkappa(m;\psi))=\bdelta_{\bkappa(m;\psi)}\bigl(\res_d\bigr)
=\psi\bigl(\res_d\bigr)=
$$
$$
=\bdelta_{\bkappa(n;\psi)}\bigl(\res_d\bigr)
=\res_d(\bkappa(n;\psi))\,,
$$
as was to be proved.
\end{proof}
\section{The ring of characters}

\subsection {Convolutions and reflection in $\ms G$} Denote by
 $\sym(\ms A\!\otimes\ms A)$ the algebra consisting of all
periodic (with respect to both arguments) functions $w:\Z^2\to\C$
such that $w(x,y)=w(y,x)$. If a function $w$ is $p$--periodic with
respect to $x$, it will be $p$--periodic with respect to $y$ as well
(and conversely). Indeed, $w(x,y+p)=w(y+p,x)=w(y,x)=w(x,y)$. Let us
agree to denote the value $\phi(u)$ of a functional $\phi\in\mc A^*$
by $\phi_xu(x)$. Thus, the variable $x$ twice repeated in this
expression is \emph{umbral} i.\,e. the letter $x$ can be replaced by
any other one. Fixing the value of one of the arguments one can
consider the function $w$ as some element of $\ms A$.
\begin{defn}
Define the \emph{direct product} $\phi\otimes\psi$ of two
functionals $\phi$ and $\psi$ by the equality
$(\phi\otimes\psi)w=\phi_x\psi_yw(x, y)$.
\end{defn}
This is a linear functional on the algebra
$\sym(\ms A\!\otimes\ms A)$.

\begin{prop}\label{commut}
The direct product is commutative.
\end{prop}

\begin{proof}
Let $w$ be an arbitrary $p$--periodic function taken from the
algebra $\sym(\ms A\!\otimes\ms A)$. Twice applying the equality
$(3)$ from Proposition\,\ref{indicators}, one can represent this
function as
\begin{equation}\label{oops}
w(x,y)=\sum_{j=0}^{p-1}w(j,y)e_j^p(x)=
\sum_{j=0}^{p-1}\sum_{k=0}^{p-1}w(j,k)e_j^p(x)e_k^p(y)\,.
\end{equation}
Let us apply the linear functionals $\phi\otimes\psi$ and
$\psi\otimes\phi$ to both sides of this equality:
$$
(\phi\otimes\psi)w=\sum_{j=0}^{p-1}\sum_{k=0}^{p-1}
w(j,k)\phi(e_j^p)\psi(e_k^p)\,.
$$
$$
(\psi\otimes\phi)w=\sum_{j=0}^{p-1}\sum_{k=0}^{p-1}
w(j,k)\psi(e_j^p)\phi(e_k^p)\,.
$$
By changing the index names in the second equality and using
the symmetry of the function $w$, we find that
\begin{multline*}
(\psi\otimes\phi)w=\sum_{k=0}^{p-1}\sum_{j=0}^{p-1}
w(k,j)\psi(e_k^p)\phi(e_j^p)=\\
=\sum_{j=0}^{p-1}\sum_{k=0}^{p-1}w(k,j)\phi(e_j^p)\psi(e_k^p)=
\sum_{j=0}^{p-1}\sum_{k=0}^{p-1}w(j,k)\phi(e_j^p)\psi(e_k^p)\,.
\end{multline*}
Hence, $\psi\otimes\phi=\phi\otimes\psi$.
\end{proof}
\noindent We shall define operators $s$ and $p$ from $\ms A$ into
$\sym(\ms A\!\otimes\ms A)$ and an operator $n$ from $\ms A$ into
$\ms A$. For any $u\in\ms A$ we set
$$
(su)(x,y)=u(x+y)\,;\quad(pu)(x,y)=u(xy)\,;\quad(nu)(x)=u(-x)\,.
$$
It is clear that these operators are multiplicative:
$$
s(uv)=(su)(sv)\,;\quad p(uv)=(pu)(pv)\,;\quad n(uv)=(nu)(nv)\,.
$$
\begin{defn}
We call the functionals $\phi\oplus\psi=(\phi\otimes\psi)\,s$ and,
respectively, $\phi\odot\psi=(\phi\otimes\psi)\,p$ a \emph{plus--}
and a \emph{dot--convolutions} of $\phi$ and $\psi$. The functional
$\ominus\phi=\phi\,n$ is called a \emph{reflection} of $\phi$.
\end{defn}

\subsection{Ring properties of the operations in $\ms G$}

\begin{prop}\label{associate}
Both convolutions in the previous definition are
commutative and associative.
\end{prop}
\begin{proof}
The commutativity of operations is an obvious corollary from
Proposition\,\ref{commut}. The associativity of the
dot--convolution follows from the equalities
$$
\begin{aligned}
(\phi\odot(\chi\odot\psi))(u)&=\phi_x(\chi\odot\psi)_s\,u(xs)=
\phi_x\chi_y\psi_z\,u(x(yz))\\
((\phi\odot\chi)\odot\psi)(u)&=(\phi\odot\chi)_s\psi_z\,u(sz)=
\phi_x\chi_y\psi_z\,u((xy)z)
\end{aligned}
$$
and from the associativity of the multiplication in $\Z$.
The associativity of the plus--convolution can be established quite
similarly .
\end{proof}

\begin{prop}\label{neutral}
The functionals $\theta=\bdelta_0$ and $\varepsilon=\bdelta_1$ are
neutral elements of the operations $\oplus$ and $\odot$
respectively.
\end{prop}
\begin{proof}
$$
\begin{aligned}
(\phi\oplus\theta)(u)&=\phi_x\theta_yu(x+y)=
\phi_xu(x+0)=\phi_xu(x)=\phi(u)\,,\\
(\phi\odot\varepsilon)(u)&=\phi_x\varepsilon_yu(xy)=
\phi_xu(x\cdot1)=\phi_xu(x)=\phi(u)\,,
\end{aligned}
$$
i.\,e. $\theta$ and $\varepsilon$ --- are neutral elements of
the corresponding convolutions.
\end{proof}
\begin{prop}\label{mult}
The direct product of characters of the algebra $\ms A$ is a
character of the algebra $\sym(\ms A\!\otimes\ms A)$.
\end{prop}

\begin{proof}
For any two functions $u$ and $v$ from the algebra $\ms A$, in view
of the multiplicativity of the functionals $\phi$ and $\psi$, we have
\begin{multline*}
(\phi\otimes\psi)(uv)=\phi_x\psi_yu(x,y)v(x,y)=
\phi_x((\psi_yu(x,y))(\psi_yv(x,y))=\\
=\phi_x\psi_yu(x,y)\phi_x\psi_yv(x,y)=
(\phi\otimes\psi)(u)(\phi\otimes\psi)(v)\,.
\end{multline*}
The multiplicativity of the functional $\phi\otimes\psi$ is thus
proved.
\end{proof}

\begin{prop}\label{stable}
The set $\ms G$ of all characters of the algebra $\ms A$ is closed
under the binary operations $\oplus$ and $\odot$, and under the
unary operation $\ominus$ of the reflection.
\end{prop}
\begin{proof}
If $\phi$ and $\psi$ are elements of $\ms G$ then, by the
multiplicativity of the operator $s$ and according to
Proposition\,\ref{mult}, we have
\begin{multline*}
(\phi\oplus\psi)(uv)=((\phi\otimes\psi)s)(uv)=
(\phi\otimes\psi)(su\,sv)=\\
=(\phi\otimes\psi)(su)(\phi\otimes\psi)(sv)=
(\phi\oplus\psi)(u)(\phi\oplus\psi)(v)\,.
\end{multline*}
Analogously, by the multiplicativity of the operator $p$, we have
\begin{multline*}
(\phi\odot\psi)(uv)=((\phi\otimes\psi)p)(uv)=
(\phi\otimes\psi)(pu\,pv)=\\
=(\phi\otimes\psi)(pu)(\phi\otimes\psi)(pv)=
(\phi\odot\psi)(u)(\phi\odot\psi)(v)\,.
\end{multline*}
The multiplicativity of the functionals $\phi\oplus\psi$ and
$\phi\odot\psi$, and the closeness of $\ms G$ under both
convolutions is thus proved. Since the operator $n$ is also
multiplicative, for any character $\phi$ we obtain
$$
(\ominus\phi)(uv)=(\phi n)(uv)=\phi(nu\,nv)=\phi(nu)\phi(nv)=
(\ominus\phi)(u)(\ominus\phi)(v).
$$
Consequently, the functional $\ominus\phi$ is a character.
\end{proof}
\begin{prop}\label{group}
The set $\ms G$ of all \,characters of the algebra $\ms A$ is a
commutative ring with the operation of addition $\oplus$,
multiplication $\odot$, zero $\theta$ and the unit element $\varepsilon$.
\end{prop}
\begin{proof}
Taking into account Propositions\,\ref{associate}, \ref{neutral},
\ref{stable} it remains to prove the existence of an opposite for
any character $\chi$ and to verify the distributivity
\begin{equation}\label{distr}
(\phi\oplus\psi)\odot\chi=(\phi\odot\chi)\oplus(\psi\odot\chi)\,.
\end{equation}
Let us show that the additive convolution of a character $\chi$
with another character $\ominus\chi$ results in an additively
neutral character $\theta$. Indeed, for an arbitrary periodic
function $u$ we have
$$
(\chi\oplus(\ominus\chi))(u)=\chi_x(\ominus\chi)_yu(x+y)=
\chi_x\chi_yu(x-y)\,.
$$
Let $p$ be the period of the function $u$. As it follows from
Proposition\,\ref{delta}, there exists an integer
$k_p=\bkappa(p;\chi)$ such that for the arbitrary element $v$ from
$\ms A_p$ the value of the functional $\chi$ is $v(k_p)$. Hence
$$
\chi_x\chi_yu(x-y)=\chi_xu(x-k_p)=
u(k_p-k_p)=u(0)=\theta(u)\,.
$$
Thus, $\chi\oplus(\ominus\chi)=\theta$. Let us verify the
distributivity for an arbitrary function $u\in\ms A_p$\,. Applying
first the left--hand side of the equality \eqref{distr} to $u$:
\begin{multline*}
((\phi\oplus\psi)\odot\chi)(u)=
(\phi\oplus\psi)_s\chi_tu(st)=
(\phi\oplus\psi)_su(s\bkappa(p;\chi))=\\
=\phi_x\psi_yu((x+y)\bkappa(p;\chi))=
u((\bkappa(p;\phi)+\bkappa(p;\psi))\bkappa(p;\chi))\,,
\end{multline*}
and then its right--hand side:
\begin{multline*}
(\phi\odot\chi)\oplus(\psi\odot\chi)(u)=
(\phi\odot\chi)_s(\psi\odot\chi)_tu(s+t)=\\
=(\phi\odot\chi)_s\psi_y\chi_zu(s+yz)=
(\phi\odot\chi)_su(s+\bkappa(p;\psi)\bkappa(p;\chi))=\\
=\phi_x\chi_zu(xz+\bkappa(p;\psi)\bkappa(p;\chi))=
u(\bkappa(p;\phi)\bkappa(p;\chi)+\bkappa(p;\psi)\bkappa(p;\chi))\,,
\end{multline*}
we see that the relation \eqref{distr} holds.
\end{proof}
\begin{defn}
In what follows, the binary operation $\phi\oplus(\ominus\psi)$
will be written as $\phi\ominus\psi$.
\end{defn}
\section{Topological ring $(\ms G,\bgamma)$.}

\subsection{Cluster topology in the ring $\ms G$.}

Let $p$ be an arbitrary positive integer.

\begin{defn}
We shall say that elements $\phi$ and $\psi$ of the ring of
characters are \emph{$p$--equivalent} if the functionals $\phi$
and $\psi$ coincide on the subalgebra $\ms A_p$. Denote by
$V_p(\psi)$ the $p$--equivalence class containing the character
$\psi$. These classes will be called \emph{clusters}.
\end{defn}

\begin{prop}\label{separate}
The character $\psi$ is the unique element which belongs to all
clusters $V_n(\psi)$.
\end{prop}

\begin{proof}
Let $\phi\in V_n(\psi)$ for all positive integers $n$. Let us show
that $\phi=\psi$. Since each function $u$ from $\ms A$ is
periodic, it belongs to a subalgebra $\ms A_n$. But as
$\phi|\ms A_n=\psi|\ms A_n$, we have $\phi(u)=\psi(u)$. As $p$ is
an arbitrary function, the functionals $\phi$ and $\psi$ coincide
on the whole algebra $\ms A$.
\end{proof}

\begin{prop}\label{basis}
For any character $\psi$ and positive integers $m$ and $n$ we have
$ V_{m\lor n}(\psi)\subseteq V_m(\psi)\cap V_n(\psi) $. where
$m\lor n$ is the least common multiple of $m$ and $n$.
\end{prop}

\begin{proof}
It is obvious that every $m$--periodic function $u$ is also
$(m\lor n)$--periodic. Therefore the coincidence of the
functionals $\phi$ and $\psi$ on the subalgebra $\ms A_{m\lor n}$
implies their coincidence on the subalgebra
$\ms A_m\subseteq\ms A_{m\lor n}$.
In other words, $V_{m\lor n}(\psi)\subseteq V_m(\psi)$.
Analogously, $V_{m\lor n}(\psi)\subseteq V_n(\psi)$.
\end{proof}

Propositions \ref{separate} and \ref{basis} mean that the system
of clusters $V_n(\psi)$ can be considered as a base of neighborhoods
of the point $\psi\in\ms G$ which defines a Hausdorff topology
on $\ms G$.

\begin{defn}
A topology defined by the base of neighborhoods $V_n(\psi)$ in the
ring $\ms G$ will be called a \emph{cluster topology}.
\end{defn}

In the following we will denote the clusters $V_n(\theta)$ by $V_n$.

\begin{prop}\label{symmetry}
Any cluster $V_n$ is a symmetric subset of $\ms G$.
\end{prop}

\begin{proof}
We have to show that if $\phi$ belongs to $V_n$ then $\ominus\phi$
belongs to $V_n$ as well. The cluster $V_n$ consists of the
elements of the ring $\ms G$ which act as the zero element on
$n$--periodic functions: $\phi(u)=\theta(u)=u(0)$. Hence for these
functions we have $(\ominus\phi)(u)=\phi_xu(-x)=u(0)$ and
$\ominus\phi\in V_n$.
\end{proof}

\begin{prop}\label{add}
Each cluster $V_n$ is a subgroup in the additive group of the ring
$\ms G$.
\end{prop}

\begin{proof}
In view of Proposition \,\ref{symmetry}, it is sufficient to show
that the cluster $V_n$ is closed under additive convolution. Let
$\phi$ and $\psi$ be elements of $V_n$. For any $n$--periodic
function $u$ we have
$$
(\phi\oplus\psi)(u)=\phi_x\psi_yu(x+y)=\phi_xu(x+0)=u(0)\,,
$$
so that $\phi\oplus\psi\in V_n$.
\end{proof}

\begin{prop}\label{ideal}
The cluster $V_n$ is an ideal of the ring $\ms G$.
\end{prop}

\begin{proof}
Let $\phi$ be an arbitrary character, and let $\psi$ be an element
of the cluster $V_n$. For each $n$--periodic function $u$ we have
$$
(\phi\odot\psi)(u)=\phi_x\psi_yu(xy)=\phi_xu(0)e(x)=
u(0)\phi(e)=u(0)\,,
$$
so $\phi\odot\psi$ is an element of $V_n$.
\end{proof}

\begin{prop}\label{translations}
$V_n(\psi)$ is a $\psi$--translation of the cluster $V_n$.
\end{prop}

\begin{proof}
Recall that the $\psi$--translation $(V_n)_\psi$ of the subset
$V_n\subseteq\ms G$ is defined by the equality
$(V_n)_\psi=\{\psi\oplus\eta:\eta\in V_n\}$.
Hence $\phi\in(V_n)_\psi$ if and only if
$(\phi\ominus\psi)|\ms A_n=\theta|\ms A_n$.
Since $V_n$ is symmetric, without loss of generality we can
consider
$
0\le\bkappa(n;\psi)\le \bkappa(n;\phi)\le n-1
$.
The fact that $\phi$ belongs to the set $(V_n)_\psi$ is equivalent
to the fact that for any function $u$ from $\ms A_n$
$$
u(\bkappa(n;\phi)-\bkappa(n;\psi))=\phi_x\psi_yu(x-y)=
(\phi\ominus\psi)(u)=\theta(u)=u(0)\,.
$$
Since the subalgebra $\ms A_n$ separates the points
$0,\ldots,n-1$, the values $\bkappa(n;\phi)$ and $\bkappa(n;\psi)$
coincide, which is equivalent to the coincidence of the characters
$\phi$ and $\psi$ on the subalgebra $\ms A_n$, i.e., to the
condition $\phi\in V_n(\psi)$.
\end{proof}

\begin{prop}\label{K-agree}
The cluster topology is compatible with the ring structure of
$\ms G$.
\end{prop}
\begin{proof}
As noted in Proposition \,\ref{translations}, the base of
neighborhoods of an arbitrary point $\psi$ from the ring $\ms G$
consists of the zero neighborhoods translations which, by Proposition
\,\ref{symmetry}, are symmetric subsets of the ring $\ms G$. It
remains to prove that the clusters $V_n$ satisfy the following
conditions from Section {\bf\ref{base}}:
\begin{enumerate}
  \item The zero of the ring $\ms G$ is the unique element
        which belongs to all $V_n$.
  \item The intersection $V_m\cap V_n$ includes some neighborhood
        of zero.
  \item For each neighborhood $V_n$ there exists a neighborhood $V_m$
        such that $\{\phi\oplus\psi:\phi,\psi\in V_m\}\subseteq V_n$.
  \item For each neighborhood $V_n$ and for any $\psi\in\ms G$ one
        can find a neighborhood $V_m$ such that
        $\{\phi\odot\psi:\phi\in V_m\}\subseteq V_n$.
  \item For any neighborhood $V_n$ there exists a neighborhood $V_m$
        such that
        $\{\phi\odot\psi:\phi\odot\psi\in V_m\}\subseteq V_n$.
\end{enumerate}
It follows from Proposition\,\ref{separate} that the character
$\theta$ is the only element belonging to all clusters $V_n$,
while Proposition\,\ref{basis} implies the condition \,$(2)$:
$V_{m\lor n}\subseteq V_m\cap V_n$. Next, in view of
Proposition\,\ref{add}, $V_n$ is a subgroup in the additive group
of the ring $\ms G$, therefore the condition\,$(3)$ is satisfied,
for instance, if $m=n$. Finally, $(4)$ and $(5)$, for $m=n$,
follow from the fact that, in view of Proposition\,\ref{ideal},
$V_n$ is an ideal of the commutative ring $\ms G$.
\end{proof}

\subsection{Cluster topology as the Gelfand topology.}

Let us correlate the cluster topology and the Gelfand topology
$\bgamma$ which, according to \eqref{neighborhood}, can be
determined by the subbase of the neighborhoods of the points
$\psi\in\ms G$:
$$
U_R(u;\psi)=\{\phi\in\ms G:\,|\phi(u)-\psi(u)|<R\}
\quad(u\in\ms A,\, R>0)\,.
$$
\begin{prop}
The cluster $V_n(\psi)$ can be represented in the form
\begin{equation}\label{representation}
V_n(\psi)=\bigcap_{j=0}^{n-1}U_R(e^n_j;\psi)\,,
\end{equation}
where $R$ is an arbitrary number from the interval $(0;1)$, and
the functions $e_j^n$ are indicators of the rasters $R_n(j)$
considered in Proposition\,\rm\ref{indicators}.
\end{prop}

\begin{proof}
First we show that $V_n(\psi)$ belongs to every neighborhood
$U_{\!R}(e^n_j;\psi)$ for any positive $R$. Since $e^n_j$ are
elements of the subalgebra $\ms A_n$, where the functionals from
$V_n(\psi)$ coincide, then for any $\phi\in V_n(\psi)$ we have
$$
\phi\in\{\phi:|\phi(e^n_j)-\psi(e^n_j)|=0\}\subseteq
\{\phi:|\phi(e^n_j)-\psi(e^n_j)|<R\}=U_{\!R}(e^n_j;\psi)\,.
$$
Thus, the right--hand side of the equality \eqref{representation}
incudes the left--hand one. Now we prove the converse inclusion.
Suppose that character $\phi$ belongs to all the neighborhoods
$U_{\!R}(e^n_j;\psi)$. According to Proposition\,\ref{delta}, its
contraction to the subalgebra $\ms A_n$ has the form
$\phi(u)=u(\bkappa(n;\phi))$, where $0\le\bkappa(n;\phi)\le n-1$.
Consequently,
$$
\phi(e^n_j)=e^n_j(\bkappa(n;\phi))=\delta_{j\bkappa(n;\phi)}\,.
$$
Since $\psi$ belongs to every neighborhood $U_{\!R}(e^n_j;\psi)$, then
$\psi(e^n_j)=\delta_{j\bkappa(n;\psi)}$. Thus, for all $0\le j\le
n-1$ the following inequalities hold:
$$
|\delta_{j\bkappa(n;\phi)}-\delta_{j\bkappa(n;\psi)}|=
|\phi(e^n_j)-\psi(e^n_j)|<R\le1\,,
$$
which mean that for any $j$ the integer
$\delta_{j\bkappa(n;\phi)}-\delta_{j\bkappa(n;\psi)}$ can only be
zero. It follows that $\bkappa(n;\phi)=\bkappa(n;\psi)$, and
therefore
$$
\phi(u)=u(\bkappa(n;\phi))=u(\bkappa(n;\psi))=\psi(u)
$$
on all $u$ from $\ms A_n$ and character $\phi$ belongs to the
cluster $V_n(\psi)$.
\end{proof}

\begin{prop}\label{G-agree}
The Gelfand topology $\bgamma$ in the ring of characters $\ms G$
coincides with the cluster topology and thus is consistent with
the ring structure of $\ms G$.
\end{prop}

\begin{proof}
In the Gelfand topology $\bgamma$, finite intersections of the
type $U=\bigcap U_{R_j}(u_j;\psi)$ form a basis of neighborhoods of
the point $\psi$, where $u_j$ are arbitrary periodic functions on
$\Z$, and $R_j$ are arbitrary positive numbers. The
equality~\eqref{representation} shows that each cluster
$V_n(\psi)$ belongs to the base of neighborhoods of the point
$\psi$. Consequently, the cluster topology is weaker than the
topology $\bgamma$. Conversely, let $p$ be the common period of
the functions $u_j$ which form the neighborhood $U$ of the element
$\psi$. For $\phi$ from the cluster $V_p$ we have
$$
\phi\in\{\phi:|\phi(u_j)-\psi(u_j)|=0\}\subseteq
\{\phi:|\phi(u_j)-\psi(u_j)|<R_j\}=U_{\!R_j}(u_j;\psi)
$$
for all values of the index $j$, so that $V_p\subseteq U$, and the
topology $\bgamma$ is weaker than the cluster topology.
\end{proof}
\section{Introduction into polyadic analysis}

\subsection{Ring $\ms P$ of polyadic numbers}

Let $\ms C$ be the commutative ring of all sequences
$\alpha:\N\to\Z$ with pointwise operations of addition and
multiplication induced from $\Z$.

\emph{A constant} $\alpha$ coincides with the result of its
\emph{shift $S$} defined by the equality
$(S\alpha)_n=\alpha_{n+1}$. Associating to every integer $m$ the
constant of the same value, we can consider the ring $\Z$ as a
subring of all constants from the ring $\ms C$. Without risk of
confusion, we shall use the same notation for constants and their
values.

As usually, a statement is said to be true \emph{for almost all}
positive integers $n$ if it is false only for a finite number
of $n\in\N$.
\par
Let us introduce a congruence relation in the ring of sequences
$\ms C$.
\begin{defn}
We shall write $\alpha\equiv\beta$ for sequences $\alpha$ and
$\beta$ if for each positive integer $n$ the congruences
$\alpha_k\equiv\beta_k\mod n$ hold for almost all $k\in\N$.
\end{defn}
In particular, this definition implies that any two sequences with
a finite number of different elements are congruent, and a sequence
is congruent with the null constant if for each positive integer
$n$ almost all its elements are divisible by $n$. This is
is an obvious analogue of the fact that zero is divisible by any
positive integer. Thus a sequence congruent with the null constant
will be called a \emph{$0$--sequence}. It is easy to see that
$0$--sequences form an ideal $\ms C_0$ of the ring $\ms C$ stable
under the shift $S$ which makes it possible to apply the shift
operation $S$ to elements of the quotient--ring $\ms C/\ms C_0$.
In addition, the elements $\ms C_0$ of the form $\alpha p$
$(p\in\N)$ can be divided by $p$: if almost all components of the
sequence $\alpha p$ are divisible by an arbitrary preassigned
$np$, then almost all components of the sequence $\alpha$ are
divisible by an arbitrary preassigned $n$.

\begin{defn}
By \emph{polyadic number} we will call any constant from the
commutative ring $\ms C/\ms C_0$, i.\,e. any class $\balpha$ of
congruent sequences $\alpha$ such that $\alpha\equiv S\alpha$.
We will say that the sequence $\alpha\in\balpha$ \emph{represents}
the polyadic number $\alpha$.
\end{defn}
\noindent
Polyadic numbers form a ring denoted by $\ms P$. The ring $\ms P$
can be considered as a module over the ring $\Z$ if we define the
product $m\balpha$ of a positive integer $m$ by an arbitrary
$\balpha$ from $\ms P$ as a class of all sequences congruent to
sequences of the form $m\alpha$ ($\alpha\in\balpha$).

\begin{defn}
By associating each positive integer $m$ with the element $\mc
Z(m)$ of the ring $\ms P$ containing the constant sequence with
the value $m$, we obtain the \emph{canonical embedding} $\mc
Z:\Z\to\ms P$.
\end{defn}
This embedding whose image will be denoted by $\ms Z$ is a strict
homomorphism of the rings.
\begin{defn}
We will say that an element $\mc Z(m)$ of the ring $\ms Z$ is an
\emph{integer polyadic number with the value $m$}.
\end{defn}
\begin{prop}\label{division}
If a sequence $\beta=n\alpha$ $(n\in\N)$ represents a polyadic
number $\bbeta$, then the sequence $\alpha$ represents a polyadic
number as well.
\end{prop}
\begin{proof}
Since $n(\alpha-S\alpha)=\beta-S\beta$ is the $0$--sequence and
since in $\ms C_0$ the division by $n\in\N$ is admissible,
then $\alpha-S\alpha$ is the $0$-- sequence.
\end{proof}

\subsection{Topological ring $(\ms P,\sigma)$}.

Let us associate each positive integer $n$ with the principal ideal
$\ms G^{(n)}$ of the ring $\ms P$ defined by the following
equality
$$
\ms G^{(n)}=\{n\balpha:\balpha\in\ms P\}.
$$

\begin{prop}\label{zero}
The zero of the ring $\ms P$ is the unique element belonging to all
ideals $\ms G^{(n)}$.
\end{prop}

\begin{proof}
Let the class $\bbeta$ be common to all $\ms G^{(n)}$, and let $p$
be an arbitrary positive integer. Then, in particular,
$\bbeta\in\ms G^{(p)}$. All components of a sequence $\beta=\alpha
p$ from the class $\bbeta$ are divisible by $p$. Hence almost all
components of any sequence from the class $\bbeta$ are divisible
by any arbitrary positive integer $p$, and thus $\bbeta$ is zero
of the ring $\ms P$.
\end{proof}

\begin{prop}\label{inclusion}
If a positive integer $p$ divides a positive integer $n$ then
$\ms G^{(n)}\subseteq\ms G^{(p)}$.
\end{prop}

\begin{proof}
Let $\bgamma\in\ms G^{(n)}$, i.\,e. $\bgamma=n\bbeta$, where
$\bbeta$ is an element of the ring $\ms P$, and let $n=pd$. All
elements of the sequence $\gamma=n\beta=p(d\beta)$ from the class
$\bgamma$ are divisible by $p$. Therefore, almost all elements of
any sequence from the class $\bgamma$ are divisible by $p$, and
$\bgamma\in\ms G^{(p)}$.
\end{proof}

\begin{prop}\label{NOK}
For arbitrary positive integers $m$ and $n$ the following equality
$\ms G^{(m)}\cap \ms G^{(n)}=\ms G^{(m\lor n)}$ holds, where
$m\lor n$ is the least common multiple of the numbers $m$ and $n$.
\end{prop}

\begin{proof}
As $m$ and $n$ divide $m\lor n$, the inclusion
$\ms G^{(m\lor n)}\!\subseteq\ms G^{(m)}\cap\ms G^{(n)}$ follows
from Proposition\,\ref{inclusion}.

Let us verify the converse inclusion. Let $\bgamma\in\ms
G^{(m)}\cap\ms G^{(n)}$, then $\bgamma=m\balpha$, and
$\bgamma=n\bbeta$, where $\balpha$, $\bbeta$ are elements of the
ring $\ms P$. Denote by $d=m\land n$ the greatest common divisor
of the numbers $m$ and $n$. The following equalities hold:
\begin{equation}\label{lor}
m=dp,\  n=dq\ (\text{where}\ p\land q=1),\ m\lor n=pdq=qm=pn\,.
\end{equation}
If the elements $\balpha$ and $\bbeta$ are represented by sequences
$\alpha$ and $\beta$ respectively, then $\alpha dp\equiv\beta dq$.
In view of Proposition\,\ref{division}, $\alpha p\equiv\beta q$,
so that the difference $\alpha_kp-\beta_kq$ is divisible by
$pq$ for almost all $k$. As $p\land q=1$, for almost all $k$ there
are numbers $\alpha'_k$ and $\beta'_k$ such that
$\alpha_k=\alpha'_kq$, and $\beta_k=\beta'_kp$. Arbitrarily
extending the functions $k\mapsto\alpha'_k$ and $k\mapsto\beta'_k$
to the sequences $\alpha'$ and $\beta'$, we obtain (taking into
account equalities \eqref{lor}) the representations of the class
$\bgamma$ by the sequences $\alpha'qm=\alpha'(m\lor n)$ and
$\beta'pn=\beta'(m\lor n)$
\end{proof}
We shall introduce a Hausdorff topology in the ring of polyadic numbers
making it a topological ring.

\begin{prop}\label{Felix}
The ideals $\ms G^{(n)}$ form a (countable) base of zero
neighborhoods, defining in the ring $\ms P$ a Hausdorff topology
$\bsigma$ compatible with its algebraic structure.
\end{prop}

\begin{proof}
A base of zero neighborhoods is defined in Section\,{\bf\ref{base}}.
Condition (1) holds by virtue of the Proposition\,\ref{zero}.
Condition (2) is provided by Proposition\,\ref{lor}.
The remaining conditions hold because $\ms
G^{(n)}$ are ideals: here one and the same set $\ms G^{(n)}$ can
be used for both sets $U$ and $V$.
\end{proof}

\section{Compactness of the ring $(\ms P,\bsigma)$}

\subsection{Grids. Finite partitions of the space $(\ms P,\bsigma)$}

\begin{defn} We call by a \emph{grid} $\ms G(n;\balpha)$ with a
\emph{width} $n\in\N$ and a \emph{center} $\balpha\in\ms P$ any
$\balpha$--translation of  the ideal $\ms G^{(n)}$, i.\,e. each
subset of the ring $\ms P$ having the form
$\ms G(n;\balpha)=\{\balpha+\bgamma:\bgamma\in\ms G^{(n)}\}$.
\end{defn}
For each point $\balpha\in\ms P$ the grids $\ms G(n;\balpha)$
form a base of its neighborhoods. Let us note some elementary
properties of grids.

\begin{prop}\label{exchange}
If $\bbeta\in\ms G(n;\balpha)$, then
$\ms G(n;\balpha)=\ms G(n;\bbeta)$.
\end{prop}
\begin{proof}
Let $\bbeta=\balpha+\bgamma$, where $\bgamma\in\ms G^{(n)}$.
Let us prove the inclusion
$\ms G(n;\balpha)\subseteq\ms G(n;\bbeta)$.
If $\bxi=\balpha+\bita$, where $\bita$ is an element of the ideal
$\ms G^{(n)}$ then $\bxi=\bbeta+(\bita-\bgamma)$, where the element
$\bita-\bgamma$ also belongs to $\ms G^{(n)}$.
Therefore, $\bxi\in\ms G(n;\bbeta)$. The inverse inclusion can
be proved analogously.
\end{proof}

\begin{prop}\label{open}
Each grid $\ms G(n;\balpha)$ is open. In particular, the ideals
$\ms G^{(n)}=\ms G(n;\theta)$ are open in the topology $\bsigma$.
\end{prop}

\begin{proof}
By Proposition\,\ref{exchange}, the grid $\ms G(n;\balpha)$ is a
neighborhood $\ms G(n;\bbeta)$ for any point $\bbeta$ of this grid.
\end{proof}

\begin{prop}\label{intersection}
If $\bdelta$ belongs to both grids $\ms G(m;\balpha)$ and
$\ms G(n;\bbeta)$, then\quad
$\ms G(m;\balpha)\cap\ms G(n;\bbeta)=\ms G(m\lor n;\bdelta)$.
\end{prop}
\begin{proof}
In view of Proposition\,~\ref{exchange}, we should prove that
$$
\ms G(m;\bdelta)\cap\ms G(n;\bdelta)=\ms G(m\lor n;\bdelta)\,,
$$
i.\,e. that $\bdelta$ are translations of the ideal $\ms G^{(m\lor n)}$,
and the intersections $\ms G^{(m)}\cap\ms G^{(n)}$ coincide.
The latter statement immediately follows from
Proposition\,\ref{NOK}.
\end{proof}
\noindent
Thus, two arbitrary grids either do not intersect or their
intersection is a grid as well.

\begin{prop}\label{alternative}
If $\bbeta\in\ms G(p;\balpha)$ then the grid $\ms G(pq;\bbeta)$
for any positive integer $q$ is absorbed by the grid
$\ms G(p;\balpha)$. Otherwise, the intersection of these grids is
empty.
\end{prop}
\begin{proof}
Let $\bbeta\in\ms G(p;\balpha)$. Since $\bbeta\in\ms G(pq;\bbeta)$,
then by Proposition\,\ref{intersection}, \,
$\ms G(p;\balpha)\cap\ms G(pq;\bbeta)=\ms G(p\lor pq;\bbeta)
=\ms G(pq;\bbeta)$, therefore
$\ms G(pq;\bbeta)\subseteq\ms G(p;\balpha)$.
\par
Now let $\bbeta\notin\ms G(p;\balpha)$ and let $\bdelta$ be an
element belonging to both grids. Then we have
$(\bdelta-\balpha)\in\ms G^{(p)}$ and
$(\bdelta-\bbeta)\in\ms G^{(pq)}\subseteq\ms G^{(p)}$. Consequently,
$\bbeta-\balpha=(\bdelta-\balpha)-(\bdelta-\bbeta)\in\ms G^{(p)}$,
i.\,e. $\bbeta\in\ms G(p;\balpha)$, which contradicts our
supposition of the existence of the element $\bdelta$.
\end{proof}
\begin{prop}\label{dense}
The subring of integer polyadic numbers $\ms Z$ is dense in the
topological ring $(\ms P,\bsigma)$.
\end{prop}
\begin{proof}
Let us verify that each grid $\ms G(n;\balpha)$ contains some
element from the ring $\ms Z$. Let the polyadic number $\balpha$
be represented by a sequence of integers $\alpha$. Then there
exists $m$ such that $\alpha_{k+1}\equiv\alpha_k\mod n$ for all
$k\ge m$. By the transitivity of the congruence,
$\alpha_k\equiv\alpha_m\mod n$ for all $k\ge m$. Let $\bzeta$ be
an integer polyadic number represented by a constant sequence
$\zeta$ with the value $\zeta_k=\alpha_m$. For $k\ge m$,
$$
(\zeta-\alpha)_k=\zeta_k-\alpha_k=\alpha_m-\alpha_k\equiv0\mod n\,,
$$
hence $\bzeta-\balpha$ belongs to the grid $\ms G^{(n)}$,
which means that the integer polyadic number $\bzeta$ belongs to
the grid $\ms G(n;\balpha)$.
\end{proof}

In the ring $\ms P$ the theorem of \emph{the division with a
remainder} holds.

\begin{prop}\label{ostatok}
For a given $n\in\N$ each polyadic number $\balpha$
can be uniquely represented in the form
$\balpha=n\bgamma(\balpha)+\brho(\balpha)$,
where $\bgamma(\balpha)\in\ms P$, and $\brho(\balpha)$ is an integer
polyadic number with the value $0\le r(\balpha)\le n-1$.
\end{prop}

\begin{proof}
Let us prove the existence of this representation. In view of
Proposition\,\ref{dense}, the grid $\ms G(n;\balpha)$ contains an
element $\bzeta=\mc Z(m)$. Let $m=kn+r$, where $0\le r\le n-1$.
Then we have
$$
\bzeta=\mc Z(nk+r)=n\mc Z(k)+\mc Z(r)\,.
$$
On the other hand, the fact that the element $\bzeta$ belongs to
the grid $\ms G(n;\balpha)$ implies $\bzeta=\balpha+n\bbeta$,
where $\bbeta\in\ms P$. Comparing these two representations of
the element $\bzeta$, we conclude that $\balpha=n\bgamma+\brho$,
where $\bgamma=\mc Z(k)-\bbeta$, and $\brho=\mc Z(r)$. Now let us
prove the uniqueness. Suppose that $\balpha=n\bgamma+\mc Z(r)$,
and $\balpha=n\bgamma'+\mc Z(r')$, where  $0\le r\le r'\le n-1$.
Subtracting term by term one representation of the polyadic number
$\balpha$ from the other, we see that $\mc
Z(r'-r)=n(\bgamma-\bgamma')$, i.\,e. non--negative constant
sequence with the value $r'-r<n$ and some integer sequence
$n(\gamma_k-\gamma_k')$ represent the same polyadic number. This
is only possible when both sequences are $0$--sequences.
\end{proof}

\begin{defn}
The integer non--negative number $r(\balpha)$ defined in
Proposition\,\ref{ostatok} will be called a \emph{remainder} of
the division of a polyadic number $\balpha$ by a positive integer
$n$ and denoted by
$$
r(\balpha)=\res_n(\balpha)\,.
$$
\end{defn}

\begin{prop}

For each positive integer $N$ one can split the ring $\ms P$
into $N$ disjoint grids with width $N$:
\begin{equation}\label{partition}
\ms P=\bigcup_{r=0}^{N-1}\ms G(N;\mc Z(r))\,.
\end{equation}
\end{prop}

\begin{proof}
According to Proposition\,\ref{ostatok}, an arbitrary polyadic
number $\balpha$ can be represented in the form
$\balpha=\bbeta+\mc Z(r)$, where the number $\bbeta$ belongs to the
ideal $\ms G^{(N)}$ and $0\le r<N$. In other words, $\balpha$
belongs to $\ms G(N;\mc Z(r))$. Since this representation is unique,
the grids corresponding to different values of $r$ do not intersect.
\end{proof}

\begin{cor}
For any positive integer $n$ and $0\le k<n!$, the grid
$\ms G(n!;\mc Z(k))$ split into disjoint grids with the
width $(n+1)!$\emph{:}
\begin{equation}\label{superpart}
\ms G(n!;\mc Z(k))=\bigcup_{m=0}^n\ms G((n+1)!;\mc Z(k+n!\cdot m))\,.
\end{equation}
\end{cor}

\begin{proof}
Let us set $N=(n+1)!$ in the equality \eqref{partition}:
$$
\ms P=\bigcup_{0\le r<(n+1)!}\ms G((n+1)!;\mc Z(r))\,.
$$
Compare the intersections of the grid $\ms G(n!;\mc Z(k))$ with the
left-- and right--hand sides of this new equality:
$$
\ms G(n!;\mc Z(k))=\bigcup_{0\le r<(n+1)!}\ms G(n!;\mc Z(k))
\cap\ms G((n+1)!;\mc Z(r))\,.
$$
Here again the right-hand side is a sum of disjoint components.
Setting $p=n!$, $q=n+1$, $\balpha=\mc Z(k)$, $\bbeta=\mc Z(r)$ in
Proposition\,\ref{alternative}, we conclude that these components
are nonempty only when $\mc Z(r-k)\in\ms G^{(n!)}$, and in this
case they coincide with the grids $\ms G((n+1)!;\mc Z(r))$. It is
easy to see that all numbers $r$ from the interval $0\le r<(n+1)!$
form a progression $r=k+n!\cdot m$, where \,$0\le m\le n$.
\end{proof}

\begin{prop}
Each grid $\ms G(n;\balpha)$ is closed. In particular, the ideals
$\ms G^{(n)}=\ms G(n;\btheta)$ are closed in the topology $\bsigma$.
\end{prop}

\begin{proof}
According to the equality \eqref{partition}, the number $\balpha$
has to belong to some grid $\ms G(n;\mc Z(k))$. Hence, in view of
Proposition\,\ref{exchange}, $\ms G(n;\balpha)=\ms G(n;\mc Z(k))$.
Let us show that the grid $\ms G(n;\mc Z(k))$ is closed. By
Proposition\,\ref{open}, the grids $\ms G(n;\mc Z(r))$ are open,
and in view of the equality \eqref{partition}, the grid $\ms
G(n;\mc Z(k))$ is the complement of the (open) sum of the grids $\ms
G(n;\mc Z(r))$ where $r\ne k$. Therefore, the (coinciding) grids
$\ms G(n;\mc Z(k))$ and $\ms G(n;\balpha)$ are closed.
\end{proof}

\subsection
{Convergence. Hausdorff compact space $(\ms P,\bsigma)$}

Convergence in the topological ring $(\ms P,\bsigma)$ is defined in
the standard way.

\begin{defn}
A sequence of polyadic numbers $\balpha^{(n)}$ \emph{converges} to
the limit $\balpha$, if each neighborhood of zero $\ms G^{(p)}$
contains almost all differences $\balpha^{(n)}-\balpha$.
\end{defn}

It is clear, that the ideal $\ms G^{(p)}$ has to contain almost all
differences\quad
$
\balpha^{(n)}-\balpha^{(n+1)}=
(\balpha^{(n)}-\balpha)-(\balpha^{(n+1)}-\balpha)\in\ms G^{(p)}
$.
The converse is also true which follows from the proposition stated
below.

\begin{prop}\label{criterium}
A sequence of polyadic numbers $\balpha^{(n)}$ converges in
$(\ms P,\bsigma)$ if and only if
$\balpha^{(n)}-\balpha^{(n+1)}\to\btheta$.
\end{prop}

\begin{proof}
Let $\balpha^{(n)}-\balpha^{(n+1)}\to\btheta$. Let us represent
each polyadic number $\balpha^{(n)}$ by an arbitrary integer
sequence $\alpha^{(n)}$. By the definition of polyadic numbers,
for any $n$ one can find $i(n)$ such that
\begin{equation}\label{i(n)}
k\ge i(n)\Rightarrow\alpha^{(n)}_k\equiv\alpha^{(n)}_{i(n)}\mod n!\,.
\end{equation}
Without loss of generality we can consider the sequence $i(n)$ as
\emph{strictly increasing} so that $i(n)\ge n$. Let us show that
the sequence $\alpha_n=\alpha^{(n)}_{i(n)}$ represents some
polyadic number $\balpha$. For an arbitrary positive integer $p$
we verify that the difference
$\alpha^{(n)}_{i(n)}-\alpha^{(n+1)}_{i(n+1)}$ is divisible by $p$
for almost all $n$. The condition
$\balpha^{(n)}-\balpha^{(n+1)}\to\btheta$ implies the existence of
$n(p)\ge p$ such that for $n\ge n(p)$ the differences
$\balpha^{(n)}-\balpha^{(n+1)}$ belong to the neighborhood of zero
$\ms G^{(p)}$, i.\,e. for $k\ge j(n)$ all remainders
$\alpha^{(n)}_k-\alpha^{(n+1)}_k$ are divisible by $p$. Now let
$k(n)=\sup\{i(n),i(n+1),j(n)\}$ and represent the difference
$\alpha_n-\alpha_{n+1}=\alpha^{(n)}_{i(n)}-\alpha^{(n+1)}_{i(n+1)}$
in the form
$$
\left(\alpha^{(n)}_{i(n)}-\alpha^{(n)}_{k(n)}\right)+
\left(\alpha^{(n)}_{k(n)}-\alpha^{(n+1)}_{k(n)}\right)+
\left(\alpha^{(n+1)}_{k(n)}-\alpha^{(n+1)}_{i(n+1)}\right).
$$
For $n\ge n(p)\ge p$ the second summand is divisible by $p$, while
the first and the last ones are divisible by $n!$ and $(n+1)!$,
respectively. Thus they both are divisible by $p$. Therefore, the
above sequence $\alpha_n$ indeed represents some polyadic number
$\balpha$, and for $k\ge n\ge n(p)$ all differences
$\alpha_n-\alpha_k$ are divisible by $p$. Let us now verify that
$\balpha^{(n)}\to\balpha$, i.\,e. for each $n\ge n(p)$ and $k\ge
i(n)$, where the sequence $i(n)\ge n$ is defined by the condition
\eqref{i(n)}, the difference $\alpha^{(n)}_k-\alpha_k$ is
divisible by $p$. To this end, we write it in the form
$$
\left(\alpha_n-\alpha_k\right)+
\left(\alpha^{(n)}_k-\alpha^{(n)}_{i(n)}\right).
$$
As has just been noted, the first summand is divisible by $p$ and
in view of \eqref{i(n)}, the second one is divisible by $n!$. By
virtue of the inequality $n\ge n(p)\ge p$, the second summand is
also divisible by $p$.
\end{proof}

\begin{prop}
The ring $(\ms P,\bsigma)$ is a Hausdorff compact space.
\end{prop}

\begin{proof}
In Proposition\,\ref{Felix} the separability of the topology
$\bsigma$ was proved. Suppose that $(\ms P,\bsigma)$ is not a
compact. Suppose that the open sets $\ms O_\lambda$ form a covering
$\ms P$ where a finite subcovering cannot be found. Then among the
two grids of the partition
$$
\ms P=\bigcup_{r=0}^1\ms G(2!;\mc Z(r))\,.
$$
one can find a grid $\ms G(2!;\mc Z(r_2))$, which either cannot be
covered with any finite sum of sets $\ms O_\lambda$, -- otherwise, the
whole space $\ms P$ would have a finite covering of elements of
the set $\ms O_\lambda$. By Proposition\,\ref{alternative}, the grid
$\ms G(2!;\mc Z(r_2))$ can be split into $3$ disjoint grids
$$
\ms G(2!;\mc Z(r_2))=\bigcup_{r=0}^2\ms G(3!;\mc Z(r_2+2!\cdot r))\,,
$$
where at least one grid $\ms G(3!;\mc Z(r_3))$ cannot be covered
with a finite sum of sets $\ms O_\lambda$ and so on. In this way we
obtain a sequence of grids
$$
\ms G(2!;\mc Z(r_2))\supseteq\ms G(3!;\mc Z(r_3))\supseteq\cdots
\supseteq\ms G(k!;\mc Z(r_k))\supseteq\cdots\,,
$$
which cannot be covered with a finite number of elements from the
set $\ms O_\lambda$. Since $\mc Z(r_{k+1})\in\ms G((k+1)!;\mc
Z(r_{k+1}) \subseteq\ms G(k!;\mc Z(r_k))$, the difference $\mc
Z(r_{k+1})-\mc Z(r_k)$ belongs to the neighborhood of zero $\ms
G^{(k!)}$. By Proposition \ref{criterium}, the sequence of integer
polyadic numbers $\mc Z(r_k)$ converges in $(\ms P,\bsigma)$ to
some $\balpha\in\ms P$. Let $\ms O$ is an element of the covering
$\ms O_\lambda$, containing the point $\balpha$. Then $\balpha$
belongs to the open set $\ms O$ together with its neighborhood
$\ms G(n,\balpha)$. Let $m$ be such that $\mc Z(r_k)\in\ms
G(n;\balpha)$ for $k\ge m$. For $k=\max\{m,n\}$ the number $k!$ is
divisible by $n$, so $\ms G^{(k!)}\subseteq\ms G^{(n)}$ and $\ms
G(k!;\mc Z(r_k))\subseteq\ms G(n;\mc Z(r_k))$. As $\mc
Z(r_k)\in\ms G(n;\balpha)$, according to
Proposition\,\ref{exchange}, $\ms G(n;\mc Z(r_k))=\ms
G(n;\balpha)$, and the following inclusion holds: $\ms G(k!;\mc
Z(r_k))\subseteq\ms G(n;\mc Z(r_k))= \ms G(n;\balpha)\subseteq\ms
O$, which contradicts our assumption.
\end{proof}
\section{Stabilizers}

The notion of stabilizer is an interpretation of the polyadic
numbers which allows one to relate the polyadic topological ring to
the topological ring of characters.

\subsection{Basic definitions}

\begin{defn}
We shall say that the sequence \,$\alpha:\N\to\Z$\,
\emph{stabi\-lizes} a function $u:\Z\to\C$ at the
\emph{final value} $v=(u\circ\alpha)_\infty$ if for almost
all $n$ the equality $u(\alpha_n)=v$ holds.
\end{defn}
\begin{defn}
A sequence $\alpha:\N\to\Z$ is said to be a \emph{stabilizer} if
it stabilizes each function $u\in\ms A$.
\end{defn}
\noindent
It follows from these definitions that for each periodic function
$u$ one can choose an integer $n$ such that the final value
$(u\circ\alpha)_\infty$ will be $u(n)$.
A ``uniformity'' condition restricts the class of stabilizers.
\par
\begin{defn}
A stabilizer $\alpha$ is \emph{absolute} if there exists an
integer $n$ such that for all $u$ belonging to $\ms A$ the
following condition holds:
$$
(u\circ\alpha)_\infty=u(n)\,.
$$
\end{defn}

\begin{prop}
For any absolute stabilizer $\alpha$ there is a unique integer $n$
which satisfies the condition $(u\circ\alpha)_\infty=u(n)$
(uniformly with respect to $u$ from $\ms A$).
\end{prop}

\begin{proof}
Let $\alpha$ be an absolute stabilizer, and let $p$ be an
arbitrary positive integer. Let $\displaystyle u=\res_p$. By definition,
there exists an integer $n$ such that the final value
$(u\circ\alpha)_\infty$ is equal to $u(n)$. If there is another
number $n'$ which has the same property, then
$\displaystyle\res_p(n')=\res_p(n)$ and the remainder $n'-n$ is divisible by
any positive integer $p$.
\end{proof}

The above proposition allows us to introduce a function $\bnu$
which is well--defined on the set of all absolute stabilizers
$\alpha$ by the equality $(u\circ\alpha)_\infty=u(\bnu(\alpha))$.

\begin{defn}
The integer $\bnu(\alpha)$ will be called the \emph{value} of an
absolute stabilizer $\alpha$.
\end{defn}

To prove that the mapping $\bnu$ is surjective it is sufficient to
consider the constant sequence with the appropriate value.

\begin{defn}
Absolute stabilizers with zero value will be called \emph{zero
stabilizers}.
\end{defn}

Thus, a zero stabilizer $\alpha$ can be described as a sequence of
integers which stabilizes each function $u$ from $\ms A$
at the value $u(0)$. Weakening this condition by exchanging the
algebra $\ms A$ for its finite--dimensional subalgebra $\ms A_p$,
we come to the concept of a \emph{prezero} stabilizer.

\begin{defn}
A stabilizer $\alpha$ is said to be a \emph{$p$--prezero} stabilizer
if $(u\circ\alpha)_\infty=u(0)$ for any $p$--periodic function $u$.
\end{defn}

\subsection{Stabilizers and polyadic numbers}

Let us describe the concepts given in the previous section in
terms of polyadic numbers.

\begin{prop}\label{criter}
A sequence $\alpha$ is a stabilizer if and only if it represents a
polyadic number.
\end{prop}

\begin{proof}
Let $\alpha$ be a stabilizer. Then for any periodic function $u$
we have $u(\alpha_{k+1})=u(\alpha_k)$ for almost all $k$.
In particular, for an arbitrary positive integer $p$ we have
$$
\res_p(\alpha_{k+1})=\res_p(\alpha_k)\,,
$$
for almost all $k$, i.\,e. $\alpha_{k+1}\equiv\alpha_k\mod p$ for
almost all $k$, and the sequence $\alpha$ represents a polyadic
number.
\vskip 4pt\noindent
Conversely, let $\alpha\in\balpha\in\ms P$. Choosing
an arbitrary function $u\in\ms A_p$, let us
verify that the sequence $\alpha$ stabilizes $u$. Denote
$\alpha'_k=\alpha_{k+1}-\alpha_k$. As the sequence $\alpha$
represents a polyadic number, $\alpha'$ is a $0$--sequence, and
beginning from $k=n$ the equalities $\alpha'_k=pq_k$ hold. Since
$u$ is a $p$--periodic function, we have for $N\ge n$
$$
u(\alpha_N)=u\Bigl(\alpha_n+\sum_{k=n}^{N-1}\alpha'_k\Bigr)=
u\Bigl(\alpha_n+p\sum_{k=n}^{N-1}q_k\Bigr)=u(\alpha_n)\,.
$$
Consequently, $\alpha$ is a stabilizer.
\end{proof}

\begin{prop}\label{criter_0}
The class of zero stabilizers coincides with the class of
$0$--sequences representing the zero element $\btheta$ of the
commutative ring $\ms P$.
\end{prop}

\begin{proof}
Let $\alpha$ be a zero stabilizer. Then for any periodic function
$u$ almost all numbers $u(\alpha_k)$ are equal to $u(0)$. In
particular, for any $p$ the relation
$$
\res_p(\alpha_n)=\res_p(0)=0\,.
$$
holds for almost all positive integers $k$. Therefore,
$\alpha\in\ms C_0$, since $\alpha_k\equiv0$ for almost all $k$.
\vskip 4pt\noindent
Conversely, let $\alpha\in\ms C_0$, i.\,e. almost all elements of
the sequence $\alpha$ are divisible by any preassigned positive
integer. For an arbitrary $p$--periodic function $u$, let us
verify that the final value of $u$ exists
and is equal to $u(0)$. Indeed, for $k\ge n(p)$ the equalities
$\alpha_k=pq_k$ hold and in view of $p$--periodicity,
$u(\alpha_k)=u(pq_k)=u(0)$.
\end{proof}

Using the notion of a zero stabilizer, it is easy to describe the
class of absolute stabilizers.

\begin{prop}\label{criter_a}
A sequence is an absolute stabilizer if and only if it represents
an integer polyadic number.
\end{prop}

\begin{proof}
Let $\alpha$ be an absolute stabilizer. Let us introduce a
constant sequence $\alpha^c$ with the value $\bnu(\alpha)$, and
define the sequence $\alpha^o$ as the remainder $\alpha-\alpha^c$.
Let us verify that $\alpha^o$ is a zero stabilizer. For an
arbitrary function $u\in\ms A$, we define the periodic function $v$
by the equality $v(\cdot)=u(\,\cdot\,-\bnu(\alpha))$. Since
$\alpha$ is an absolute stabilizer with the value $\bnu(\alpha)$,
the sequence $\alpha^o$ stabilizes $u$ at the final
value $u(0)$:
$$
(u\circ\alpha^o)_\infty=(v\circ\alpha)_\infty=v(\bnu(\alpha))=
u(\bnu(\alpha)-\bnu(\alpha))=u(0)\,.
$$
Thus $\alpha^o$ is a zero stabilizer, and by
Proposition\,\ref{criter_0}, it is a $0$--sequence, which implies
that $\alpha$ is an integer polyadic number with the value
$\bnu(\alpha)$.

\vskip 4pt\noindent Conversely, suppose that $\alpha$ represents an integer
polyadic number with the value $m$, i.\,e.
$\alpha_k=m+\alpha^o_k$, where $\alpha^o$ is a $0$--sequence and
(by Proposition\,\ref{criter_0}) a zero stabilizer as well. Let us
verify that $\alpha$ is an absolute stabilizer with the value $m$.
For an arbitrary function $u$ from $\ms A$, we define the function
$v$ by the equality $v(\cdot)=u(\,\cdot+m)$. As $\alpha^o$ is a
zero stabilizer, the sequence $\alpha$ stabilizes the function $u$
at the final value
$(u\circ\alpha)_\infty=(v\circ\alpha^o)_\infty=v(0)=u(m)$.
\end{proof}

\begin{prop}\label{prezero}
A sequence is a $p$--prezero stabilizer if and only if it
represents an element of the ideal $\ms G^{(p)}$.
\end{prop}

\begin{proof}
Let $\alpha$ be a $p$--prezero stabilizer. By
Proposition\,\ref{criter}, it represents a polyadic number
$\balpha$. Let us show that $\balpha$ belongs to the ideal
$\ms G^{(p)}$. Indeed, by the definition of a $p$--prezero
stabilizer, for almost all $k$ we have
$$
\res_p(\alpha_k)=\res_p(0)=0\,,
$$
hence almost all elements of the sequence have the form
$\alpha_k=p\,\beta_k$, where $\beta$ is some sequence of integers.
According to Proposition\,\ref{division}, the sequence $\beta$
represents a polyadic number $\bbeta$. Therefore, $\alpha$
represents the polyadic number $\balpha=p\,\bbeta$, belonging to
the ideal $\ms G^{(p)}$.
\vskip 4pt\noindent
Conversely, suppose that a sequence $\alpha$ represents a polyadic
number $\balpha$ which belongs to the ideal $\ms G^{(p)}$,
i.\,e., for almost all $k$ elements of the sequence have the form
$\alpha_k=p\,\beta_k$, where $\beta$ is some sequence of integers.
By Proposition\,\ref{criter}, the sequence $\alpha$ is a
stabilizer. Let us verify that this stabilizer is $p$--prezero. In
view of $p$--periodicity, for $u\in\ms A_p$ the equalities
$u(\alpha_k)=u(p\,\beta_k)=u(0)$ hold for almost all $k$.
\end{proof}
\section
{Isomorphism of the rings $(\ms P,\bsigma)$ and $(\ms G,\bgamma)$}

\subsection{Algebraic isomorphism of the rings $\ms P$ and $\ms G$}

Let $\mc S$ be a functional ring of all stabilizers with pointwise
operations induced from $\C$. Let us associate each stabilizer
$\alpha$ from $\mc S$ with the functional $\pi\alpha:\ms A\to\C$
defined by the equality

\begin{equation}\label{pi}
(\pi\alpha)(u)=(u\circ\alpha)_\infty.
\end{equation}

\begin{prop}\label{morphism}
For any stabilizer $\alpha$ the functional $\pi\alpha$ is a
character.
\end{prop}

\begin{proof}
Let $e$ be the unit of the algebra $\ms A$, i.\,e. $e$ is a
constant function on $\Z$ with the value $1$. Since for any
stabilizer $\alpha$ the composition $e\circ\alpha$ is a constant
sequence with the value $1$, then
$(\pi\alpha)(e)=(e\circ\alpha)_\infty=1$ and the functional
$\pi\alpha$ is nonzero. For $z\in\C$ and any function $u\in\ms A$
we have
$$
(\pi\alpha)(zu)=((zu)\circ\alpha)_\infty=
(z(u\circ\alpha))_\infty=z(u\circ\alpha)_\infty=z(\pi\alpha)(u)\,,
$$
which proves the homogeneity of the functional $\pi\alpha$.
Using the distributive properties of the composition
$(u+v)\circ\alpha=u\circ\alpha+v\circ\alpha$ and
$(uv)\circ\alpha=(u\circ\alpha)(v\circ\alpha)$, it is easy to verify
its additivity
\begin{multline*}
(\pi\alpha)(u+v)=((u+v)\circ\alpha)_\infty=
(u\circ\alpha+v\circ\alpha)_\infty=\\
=(u\circ\alpha)_\infty+(v\circ\alpha)_\infty
=(\pi\alpha)(u)+(\pi\alpha)(v)
\end{multline*}
and multiplicativity
\begin{multline*}
(\pi\alpha)(uv)=((uv)\circ\alpha)_\infty=
((u\circ\alpha)(v\circ\alpha))_\infty=\\
=(u\circ\alpha)_\infty(v\circ\alpha)_\infty
=(\pi\alpha)(u)(\pi\alpha)(v)\,.
\end{multline*}
Thus, we have all the properties of characters.
\end{proof}

The above proposition allows us to introduce a mapping $\pi$
acting from the ring $\mc S$ into the ring $\ms G$ by the rule
$\pi:\alpha\mapsto\pi\alpha$.

\begin{prop}\label{homomorphism}
The mapping $\pi$ is a homomorphism of the ring $\mc S$ into the
ring $\ms G$, whose kernel $\ker\pi$ is the ideal $\ms C_0$ of
zero stabilizers.
\end{prop}

\begin{proof}
First we prove that the mapping $\pi$ is additive. It is
sufficient to verify that the values of the functionals
$\pi(\alpha+\beta)$ and $\pi\alpha\oplus\pi\beta$ are equal for
any periodic function $u$. Let $p$ be the period of the function
$u$. Let us quote the equality \eqref{oops} which holds for any
function $w$, $p$--periodic with respect to both its arguments:
\begin{equation}\label{again}
w(x,y)=\sum_{j=0}^{p-1}\sum_{k=0}^{p-1}w(j,k)e_j^p(x)e_k^p(y)\,.
\end{equation}
For $w(x,y)=u(x+y)$ this equality has the form
$$
u(x+y)=\sum_{j=0}^{p-1}\sum_{k=0}^{p-1}u(j+k)e_j^p(x)e_k^p(y)\,.
$$
Using this form, the composition $u\circ(\alpha+\beta)$ can be
rewritten as follows:
$$
u\circ(\alpha+\beta)=
\sum_{j=0}^{p-1}\sum_{k=0}^{p-1}u(j+k)
(e_j^p\circ\alpha)(e_k^p\circ\beta)\,.
$$
Passing to final values, we have, on one hand,
$$
(\pi(\alpha+\beta))(u)=
\sum_{j=0}^{p-1}\sum_{k=0}^{p-1}u(j+k)\,
\pi\alpha(e_j^p)\,\pi\beta(e_k^p)\,.
$$
On the other hand, since $\phi_x\psi_yu(x)v(y)=\phi(u)\psi(v)$
for arbitrary characters $\phi$, $\psi$ and functions $u$, $v$
from $\ms A$, we have
\begin{multline*}
(\pi\alpha\oplus\pi\beta)(u)=(\pi\alpha)_x(\pi\beta)_yu(x+y)=\\
=\sum_{j=0}^{p-1}\sum_{k=0}^{p-1}
u(j+k)(\pi\alpha)_x(\pi\beta)_ye_j^p(x)e_k^p(y)=\\
=\sum_{j=0}^{p-1}\sum_{k=0}^{p-1}
u(j+k)\,\pi\alpha(e_j^p)\,\pi\beta(e_k^p)\,.
\end{multline*}
The additivity of the mapping $\pi$ is thus proved. To prove its
multiplicativity it is sufficient to verify that for any periodic
function $u$ the values of the functionals $\pi(\alpha\beta)$ and
$\pi\alpha\odot\pi\beta$ coincide. For $w(x,y)=u(xy)$ the equality
\eqref{again} can be transformed to the form
$$
u(xy)=\sum_{j=0}^{p-1}\sum_{k=0}^{p-1}u(jk)e_j^p(x)e_k^p(y)\,,
$$
then the multiplicativity of the mapping $\pi$ can be proved in
the same way as its additivity.
\par
Thus we have ascertained that $\pi$ is a homomorphism of the ring
of stabilizers $\mc S$ into the ring of characters $\ms G$.
Belonging to the kernel $\ker\pi$ for a stabilizer $\alpha$ is
equivalent to the fact that $\pi\alpha$ is the neutral element
$\theta$ of the additive group of the ring $\ms G$, i.\,e., to the
fact that $\pi\alpha(u)=u(0)$ for any function $u$ from $\ms A$.
Since $\pi\alpha(u)=(u\circ\alpha)_\infty$, then
$$
\ker\pi=
\{\alpha\in\mc S:\forall u\in\ms A\ (u\circ\alpha)_\infty=u(0)\}=
\ms C_0\,,
$$
where $\ms C_0$ is the ideal of zero stabilizers.
\end{proof}

\begin{prop}\label{image}
The homomorphic image of $\mc S$ is the ring $\ms G$.
\end{prop}

\begin{proof}
Let $\psi$ be an arbitrary character. As noted in Section
{\bf\ref{almost}}, the algebra $\ms A$ is the inductive limit of
the increasing (with respect to inclusion) sequence of subalgebras
$\ms A_{n!}$. As $\ms G$ does not contain the zero functional, for
almost all $n$ the contractions $\psi_{n!}=\psi|\ms A_{n!}$ are
nontrivial and beginning from $n=n(\psi)$, we can define a
sequence $\alpha_n=\bkappa(n!;\psi)$ which we extend as zero for
smaller values of $n$. According to Proposition\,\ref{WOW!}, we
have
\begin{equation}\label{stab}
\alpha_{n+1}-\alpha_n=\bkappa((n+1)!;\psi)-\bkappa(n!;\psi)
\equiv0\mod n!\,.
\end{equation}
Let $u$ be an arbitrary function from $\ms A$ and let $p$ be its
period. We set $q$ large enough to satisfy the inequality $N=qp\ge
n(\psi)$. It follows from \eqref{stab} that for $n\ge N$ the
relation
$$
\alpha_n=\alpha_N+\sum_{r=N}^{n-1}(\alpha_{r+1}-\alpha_r)
\equiv\alpha_N\mod N!\,,
$$
holds where, as usual, the ``empty'' sum is said to be zero.
As $N!$ is the period of the function $u$, for $n\ge N$ we have
$$
u(\alpha_n)=u(\alpha_N)=u(\bkappa(N!);\psi)=\psi_{N!}(u)=\psi(u)\,.
$$
Thus, the sequence $\alpha$ constructed above stabilizes any
function $u$ to the final value $(u\circ\alpha)_\infty=\psi(u)$.
In other words, $\pi\alpha=\psi$.
\end{proof}
Taking into account Proposition\,\ref{homomorphism}, we obtain
\begin{cor}
The quotient mapping $\bpi=\pi/\ms C_0$ acts as an isomorphism of
the rings $\ms P$ and $\ms G$.
\end{cor}

\subsection{Isomorphism of the topological rings $(\ms P,\bsigma)$
and $(\ms G,\bgamma)$}

The following proposition shows that $\bpi$ is not only an
isomorphism of the rings $\ms P$ and $\ms G$, but an homeomorphism
of the topological spaces $(\ms P,\bsigma)$ and $(\ms G,\bgamma)$
as well, since here the zero neighborhoods of the ring
$(\ms P,\bsigma)$ passe to the zero neighborhoods of the
ring $(\ms G,\bgamma)$.

\begin{prop}\label{homeomorphism}
For each positive integer $p$ the cluster $V_p$ is the image of
the ideal $\ms G^{(p)}$ under the isomorphism $\bpi$.
\end{prop}

\begin{proof}
As stated in Proposition\,\ref{prezero}, the sequence $\alpha$
represents an element of the ideal $\ms G^{(p)}$ if and only if it
is $p$--prezero stabilizer. Hence, the proposition will be proved
if we show that the homomorphism $\pi$ surjectively maps
the set of $p$--prezero stabilizers onto the cluster $V_p$.

Let a stabilizer $\alpha$ be $p$--prezero. Then for an arbitrary
function $u$ from $\ms A_p$ \,
$(\pi\alpha)(u)=(u\circ\alpha)_\infty=u(0)=\theta(u)$. Thus
$\pi\alpha\in V_p$.

Conversely, let us suppose that a character belongs to the cluster
$V_p$, i.\,e., the functional $\psi$ has the value $v(0)$ for any
$p$--periodic function $v$. Then the stabilizer $\alpha$
constructed in the proof of Proposition\,\ref{image} is
$p$--prezero: $(v\circ\alpha)_\infty=\psi(v)=v(0)$.
\end{proof}

The following main result of this paper is an immediate consequence
of Proposition \ref{image} (Corollary) and Proposition
\ref{homeomorphism}:

\begin{theorem}
The ring $(\ms P,\bsigma)$ is algebraically and topologically
isomorphic to the ring $(\ms G,\bgamma)$.
\end{theorem}
This means that we can identify the polyadic topological ring
$(\ms P,\bsigma)$ of Pr\"ufer, Van Dantzig and Novoselov with the
topological ring of characters of the algebra $\ms A$ consisting
of all periodic functions on $\Z$.
\section*{Acknowledgements}

The authors are greatly indebted to Prof. V.~Shelkovich for
fruitful discussions. They are also grateful to Mrs. E.~Sorokina
for her essential help with this work.
\par
Besides, the second author would like to express his thanks to the
DFG for the financial support.

\end{document}